\documentclass[10 pt]{amsart}
\usepackage{amscd,amsmath,amsthm,amssymb}
\usepackage{enumerate,varioref}
\usepackage{epsfig}
\usepackage{graphicx}
\usepackage[all]{xy}
\newtheorem{thm}{Theorem}

\theoremstyle{definition}
\newtheorem{defn}[thm]{Definition}
\theoremstyle{remark}
\newtheorem{ex}[thm]{Example}
\newtheorem{rem}[thm]{Remark}
\numberwithin{equation}{subsection}

\newcommand{\Z}{\mathbb{Z}}
\newcommand{\C}{\mathbb{C}}

\begin{document}

\title{Matrix Corepresentations for $SL_q(N)$ and $SU_q(N)$}
\author{Clark Alexander}

\email{gcalex@temple.edu} \maketitle \tableofcontents

\section*{Introduction} The present paper is meant to give an algorithm for
computing matrix corepresentations of the quantum groups $SL_q(N)$
and $SU_q(N)$.  This is done by first computing the
corepresentations for $N=2$ as in cite[KS] then a simple
combinatorial re-indexing of basis elements leads one to a similar
method for computing $N>2$. The computations will be given
explicitly when possible, however the closed form of the
corepresentation is somewhat difficult to write down.  The paper
following this one will use these corepresentations to compute the
Haar functional on $SU_q(N)$ for all $N$.

\section{Notations and Preliminaries}
Much of the current literature on quantum groups uses $t_{ij}$ as
the notation for an element of a quantum matrix group.  Presently
this notation will be used for corepresentations and thus another
notation is necessary for elements within quantum groups.  To this
end, let $u_{ij} \in SU_q(N)$.  Clearly $u_{ij}\in SL_q(N)$ also,
but the existence of a $*$ product on $SL_q(N)$ will allow easy
transitioning between the two algebras. When necessary
$u_{ij}^{(N)}$ will be used if confusion is to arise as to which
space contains a
particular element.\\

For the duration of this article, let $SL_q(N)$ and $SU_q(N)$ be
used to denote the coordinate algebras usually denoted $\C_q[SL(N)],
\C_q[SL_N], \mathcal{O}_q(SL_N),$ or $\mathcal{O}(SL_q(N))$.  Here,
$q$ is a transcendental number between zero and one.

\subsection{Hopf Algebra Structure}
The algebras $SL_q(N)$ are given by $N^2$ generators labelled
$u_{ij}, i,j=1,\dots,N$ with the following relations
\[
u_{ij} u_{ik} = q u_{ik}u_{ij} , j<k, \quad u_{ij} u_{kj} = q
u_{kj}u_{ij},  i<k
\]
\[
[u_{ij},u_{kl} \rbrack = (q-q^{-1})u_{il}u_{kj},  i<k,j<l.
\]

The specific defining relation for $SL_q(N)$ is the \emph{quantum
determinant} relation given by
\begin{equation}
det_q = \sum_{\sigma\in {S_N}} (-q)^{\ell(\sigma)}
u_{1,\sigma(1)}\cdots u_{N,\sigma(N)} = \sum_{\sigma\in{S_N}}
(-q)^{\ell(\sigma)} u_{\sigma(1),1}\cdots u_{\sigma(N),N}=1.
\end{equation}

The coproduct $\Delta$ is given by
\[ \Delta(u_{ij}) = \sum_{k=1}^N u_{ik}\otimes u_{kj}.\]
The counit $\epsilon$ is given by
\[ \epsilon(u_{ij}) = \delta_{ij}.\]

\subsection{Quantum Determinants, Antipodes, $*$-Product}
To turn $SL_q(N)$ into $SU_q(N)$ one requires a $*$-product.  The
algebra $SU_q(N)$ is given by the $2N$ generators
$\{u_{ij},u^*_{ij}\}, i,j=1,\dots,N$ with the above relations and
additional relations arising from quantum ''minor" determinants.

\begin{defn}Let $\Omega_n =$ subsets
of $\{1,2,\dots,N \}$ containing $n$ elements.  Then for any
$I,J\in\Omega_n$ one writes $I=\{i_1,\dots,i_n \}$ and
$J=\{j_1\dots,j_n\}$ with $i_1<i_2<\cdots<i_n$ and
$j_1<j_2<\cdots<j_n.$  The \emph{minor determinate}
$\mathcal{D}_J^I$ is given by
\begin{equation}
\mathcal{D}_J^I = \sum_{\sigma\in
S_n}(-q)^{\ell(\sigma)}\prod_{k=1}^n u_{i_k,\sigma(j_k)} =
\sum_{\sigma\in S_n}(-q)^{\ell(\sigma)}\prod_{k=1}^n
u_{\sigma(i_k),j_k}
\end{equation}
\end{defn}

These minors have the nice properties that
\begin{equation}
\Delta(\mathcal{D}_J^I) =
\sum_{K\in\Omega_n}\mathcal{D}_K^I\otimes\mathcal{D}_J^K
\end{equation}
and
\begin{equation}
\varepsilon(\mathcal{D}_J^I)=\delta_{IJ}.
\end{equation}
When $I=J=\{1,2,\dots,N\}$ one says $det_q=\mathcal{D}_J^I.$\\
\indent The following relations are direct consequences of the above
formulae.
\begin{eqnarray}
\Delta(det_q) = det_q\otimes det_q, & \varepsilon(det_q)=1
\end{eqnarray}

\begin{defn} The generators $u^*_{ij}$ are given by the following
formula;
\begin{equation}
u^*_{ij} = D^I_J, I=\{1,\dots,\hat i,\dots,N\},J=\{1,\dots,\hat
j,\dots,N \}.
\end{equation}
\end{defn}

The anitpode $S$ is then given by
\begin{equation}
S(u_{ij}) = u^*_{ji}
\end{equation}

\subsection{q-Combinatorial Formulae}
\indent For $0\neq q\in \C$ and $a\in\Z$ define the $q$-integer
$$[a]_q = \frac{q^{a}-q^{-a}}{q-q^{-1}}.$$
One may also define the following entities
\begin{eqnarray*}
[m]! &=& [m][m-1]\cdots[1]\\
(a;q)_n &=& \prod_{j=0}^{n-1}(1-aq^j)\\
\left[ \begin{array}{c} m\\n\end{array}\right]_q &=&
\frac{(q;q)_m}{(q;q)_n(q;q)_{m-n}}\\
\left[ \begin{array}{c} m\\(i,j,k)\end{array}\right]_q &=&
\frac{(q;q)_m}{(q;q)_i(q;q)_j(q;q)_k}.
\end{eqnarray*}

The only relevant piece of information left to give is the analog of
the binomial theorem.
\begin{thm} Let $x,y$ be noncommuting variables such that $xy=qyx$
then the following formula holds.
\begin{equation}
(x+y)^n = \sum_{k=0}^n \left[\begin{array}{c}n\\k \end{array}
\right]_{q^{-1}}x^ky^{n-k}
\end{equation}
\end{thm}

One should never fail to realize that throughout these notes if one
considers the limit $q\rightarrow 1$ then one obtains the classical
situation.  For example $[a]_q = q^{a-1} + \cdots +q^{1-a}$ where in
fact there are $a$ copies of $q$ present.

\section{Corepresentations for $N=2$}

\begin{defn} Let $\mathcal{A}$ be a Hopf algebra with counit $\varepsilon$ and comultiplication $\Delta.$
Let $V$ be a complex vector space.  Then a corepresentation of
$\mathcal{A}$ on $V$ is a linear map $\varphi: V\rightarrow V\otimes
\mathcal{A}$ such that the two relations;
\[
(id\otimes \Delta)\circ \varphi =(\varphi\otimes id)\circ \varphi.,
\quad (id\otimes \epsilon)\circ \varphi=id
\]
are satisfied.
\end{defn}

Equivalently, one should require the following diagrams to commute.
\[
\xymatrix{V \ar[r]^{\varphi}\ar[d]_{\varphi} & V\otimes\mathcal{A}\ar[d]_{\varphi\otimes id} &  & V\ar[r]^{\simeq} \ar[d]_{\varphi}&V\otimes \C\ar[d]_{id} \\
          V\otimes\mathcal{A}\ar[r]^-{id\otimes\Delta}& V\otimes\mathcal{A}\otimes\mathcal{A} && V\otimes\mathcal{A}\ar[r]^-{id\otimes \epsilon} &V\otimes\C.      }
\]

 One immediately sees that when $\varphi$ is a corepresentation
of $\mathcal{A}$ on $V,$ $V$ is a right $\mathcal{A}$-comodule with
right coaction $\varphi.$\\
\indent In order to give a clear exposition of the quantum case it
is essential for one to first examine the classical case of matrix
corepresentations.  First consider the case of $SL(2,\C)$ and its
corepresentations.  Let $f\in \C[x,y]$ be a homogeneous polynomial
of degree $2\ell$ for $\ell\in \frac{1}{2}\mathbb{N}_0.$  Then one
defines the left and right actions of $SL(2,\C)$ on $f$ as follows:
letting $g = \left(\begin{array}{cc}a&b\\c&d \end{array}\right) \in
SL(2,\C)$
\begin{eqnarray}
(T_{\ell}^R(g)f)(x,y) = f(ax+cy,bx+dy)\nonumber\\
(T_{\ell}^L(g)f)(x,y) = f(ax+by,cx+dy)
\end{eqnarray}

Pictorially one should view these as left and right matrix actions
on vectors.  In the case $\ell=\frac{1}{2}$ the homogeneous
polynomials are simply $x$ and $y.$ Hence one obtains the actions
\begin{eqnarray*}
(T_{\ell}^R(g)f)(x,y) =& f((x,y)\left(\begin{array}{cc}a&b\\c&d \end{array}\right))\\
(T_{\ell}^L(g)f)(x,y) =& f(\left(\begin{array}{cc}a&b\\c&d
\end{array}\right)\left(\begin{array}{c}x\\y \end{array}\right))
\end{eqnarray*}

\begin{defn}
The matrices $T_{\ell}^R$ and $T^L_{\ell}$ are called \emph{matrix
corepresentations} of the matrix group $G$ in which $g\in G$
determine the coaction as above.
\end{defn}

In the case $\ell=\frac{1}{2}$ one determines that
\[ T_{1/2} = \left(\begin{array}{cc} a&b\\c&d\end{array}\right).\]

\begin{rem}To see that this definition is not trivial, consider a
homogeneous polynomial $f\in\C[x,y]$ of degree 2.  Carrying out the
computations one finds that
\[T_1 = \left(\begin{array}{ccc} a^2 &\sqrt{2}ab&b^2\\\sqrt{2}ac & ad+bc&
\sqrt{2}bd \\c^2&\sqrt{2}cd & d^2\end{array}\right). \]

There is a subtlety in the above computation, but that will be
explored more fully in the quantum case when choosing a basis for
$\mathcal{O}(\C^n_q)_{2\ell}$ matters.
\end{rem}

The case is nearly an exact analogy in the quantum setting, however
the vector spaces have changed and the bases require adjustments to
insure the left and right corepresentations match.

In $N$ dimensions the proper vector space over which one works are
denoted $\mathcal{O}(\C^N_q)$ or simply $\C^N_q$ given by:
\[
\mathcal{O}(\C^N_q) = \{x_i | i=1,\dots,N \hspace{.05 in}\mathrm{
and }\hspace{.05in} x_ix_j = q x_jx_i
\hspace{.05in}\mathrm{when}\hspace{.05in} i<j\}.
\]

The remainder of this section will concentrate only on the case
$N=2$.  The coactions are denoted $\varphi_R$ and $\varphi_L$ for
right and left coactions given by the formulae:
\begin{eqnarray}
\varphi_R(x_i) = \sum_{j=1}^2 x_j\otimes u_{j,i},\\
\varphi_L(x_i) = \sum_{j=1}^2 u_{i,j}\otimes x_j.\nonumber
\end{eqnarray}

Consider the case $\ell = 1/2$ then one sees
\begin{equation}
T_{1/2}=T^R_{1/2} = T^L_{1/2} = \left(\begin{array}{cc} u_{11} & u_{12}\\
u_{21} & u_{22} \end{array}\right).
\end{equation}

\begin{rem}
One needs to remember here that in the quantum case $T_{\ell}$ is
not actually a matrix, but will still be referred to as a
\emph{matrix corepresentation} of the quantum group $SL_q(2)$.
\end{rem}

\begin{defn} The matrix corepresentations of $SL_q(2)$ given by
$T_{\ell}$ are given in the form
\begin{equation}
T_{\ell} = \{  t^{\ell}_{ij}\}_{i,j = -\ell}^{\ell}
\end{equation}
\end{defn}

\begin{rem}
One might notice that the indices of $t_{ij}^{\ell}$ run from
$-\ell$ to $\ell$ which in the case of $T_{1/2}$ means
\[
T_{1/2} = \left(\begin{array}{cc} t^{1/2}_{-1/2,-1/2} &
t^{1/2}_{-1/2,1/2}\\
t^{1/2}_{1/2,-1/2} & t^{1/2}_{1/2,1/2} \end{array}\right) =
\left(\begin{array}{cc} u_{11}& u_{12}\\u_{21} & u_{22}
\end{array}\right).
\]
The issue of re-indexing corepresentation elements is the subject of
the appendix.  For now, suffice it to say that one wants
$t_{ij}^{\ell}$ to be symmetric about 0 in $i$ and $j$.
\end{rem}

\begin{defn}
The matrix corepresentation $T_{\ell}$ is called the \emph{spin
$\ell$ corepresentation} of $SL_q(2)$.
\end{defn}

\indent One begins a fit of problematic computations when one takes
the homogeneous basis of $\mathcal{O}(\C^N_q)_{2\ell}$ to simply be
the list $\{ x_1^{2\ell}, x_1^{2\ell-1}x_2, x_1^{2\ell-1}x_3,\dots,
x_N^{2\ell}\}.$ The issue here is normalization.  Consider for a
moment the case $\ell=1$ on $SL_q(2).$  Running the computation with
the faulty basis one procures the equations
\begin{eqnarray*}
T_1^R (x_1^2,x_1x_2,x_2^2) =& (x_1^2,x_1x_2,x_2^2)\otimes\left(
\begin{array}{ccc}
u_{11}^2&u_{11}u_{12}&u_{12}^2\\
(1+q^{-2})u_{11}u_{21}&u_{11}u_{22}+q^{-1}u_{12}u_{21}&(1+q^{-2})u_{12}u_{22}\\
u_{21}^2&u_{21}u_{22}&u_{22}^2
\end{array}\right) \\
T_1^L \left(\begin{array}{c}x_1^2\\x_1x_2\\x_2^2\end{array}\right)
=& \left(
\begin{array}{ccc}
u_{11}^2&(1+q^{-2})u_{11}u_{12}&u_{12}^2\\
u_{11}u_{21}&u_{11}u_{22}+q^{-1}u_{12}u_{21}&u_{12}u_{22}\\
u_{21}^2&(1+q^{-2})u_{21}u_{22}&u_{22}^2
\end{array}\right)\otimes\left(\begin{array}{c}x_1^2\\x_1x_2\\x_2^2\end{array}\right)
\end{eqnarray*}

There is an obvious problem in that $T_1^R\neq T_1^L,$ however $T_1$
is defined to be the matrix of coefficients from $T_1^R$ and
$T_1^L.$  These two corepresentations are required to match!\\
\indent One will arrive at the same problem in the classical case by
assuming the analogous faulty basis.  The solution is to renormalize
the basis in the following way: consider the binomial equation in
the commutative case
\[
(x+y)^k = \sum_{j=0}^k \left(\begin{array}{c}k\\j
\end{array}\right)x^jy^{k-j}.
\]
And more generally
\[
(x_1 + \cdots + x_n)^k =
\sum_{j_1+\cdots+j_n=k}\left(\begin{array}{c}k\\(j_1,\dots,j_n)
\end{array}\right)\prod_{i=1}^n x_{i}^{j_i}.
\]
These equations suggest that a more suitable basis for homogeneous
polynomials involves a binomial coefficient multiplier for the mixed
terms.  In fact after some delineation one will discover that a
proper basis for computing corepresentations in the commutative case
is $\{ \left(\begin{array}{c}k\\(j_1,\dots,j_n)
\end{array}\right)^{1/2}\prod_{i=1}^n x_{i}^{j_i} \}$ with a square
root hitting the binomial and multi-nomial coefficients.  Of course
this leads one to conjecture that the basis in the quantum case will
result in a similar basis with $q$-binomial coefficients, however
the elements used in computing the matrix corepresentations do not
commute with a factor of $q,$ but instead $q^{-2}.$
\begin{ex}
Consider again the case of $\ell=1$ on $SL_q(2).$
\begin{eqnarray*}
\varphi_R(x_1x_2)&=&\varphi_R(x_1)\varphi_R(x_2)\\ &=& x_1^2\otimes
u_{11}^2 + x_1x_2\otimes u_{11}u_{21} + x_2x_1\otimes u_{21}u_{11} + x_2^2\otimes u_{21}^2\\
&=&x_1^2\otimes u_{11}^2 + (1+q^{-2})xy\otimes u_{11}u_{21} +
x_2^2\otimes u_{21}^2.
\end{eqnarray*}
Here the appropriate variables to consider are the $x_i\otimes
u_{j,k}$ for some $i,j,k.$  Notice that $(x_1\otimes
u_{11})(x_2\otimes u_{21}) = q^{2}(x_2\otimes u_{21})(x_1\otimes
u_{11}).$
\end{ex}

When the smoke clears the resulting appropriate basis for
$\mathcal{O}(\C^N_q)_{2\ell}$ is
\[
\{\left(\begin{array}{c}2\ell\\(j_1,\dots,j_N)
\end{array}\right)^{1/2}_{q^{-2}} \prod_{j=1}^N x_{i}^{j_i} \}
\]
where $\sum_i j_i=2\ell.$

Running the computation again in the new basis will yield
\begin{equation}
T_1 = \left(\begin{array}{ccc} u_{11}^2 & (1+q^{-2})^{1/2}
u_{11}u_{12} & u_{12}^2\\
(1+q^{-2})^{1/2}u_{11}u_{21} & u_{11}u_{22}+q^{-1}u_{12}u_{21} &
(1+q^{-2})^{1/2} u_{12}u_{22}\\
u_{21}^2 & (1+q^{-2})^{1/2}u_{21}u_{22}& u_{22}^2
\end{array}\right).
\end{equation}

With the technology built here, an algorithm for computing matrix
corepresentations for $SL_q(2)$ becomes apparent.
\begin{enumerate}
\item Given $\ell\in \frac{1}{2}\mathbb{N}$ choose as a basis for
$\mathcal{O}(\C^2_q)_{2\ell}$ the set
\[
\{\left(\begin{array}{c}2\ell\\j
\end{array}\right)^{1/2}_{q^{-2}}  x_{1}^{j}x_2^{2\ell-j} \}.
\]\\
\item Calculate $T_{\ell}(\left(\begin{array}{c}2\ell\\j
\end{array}\right)^{1/2}_{q^{-2}}  x_{1}^{j}x_2^{2\ell-j}).$  At this point $T^R = T^L$. \\
\item In order to write down $t_{ij}^{\ell}$ explicitly one needs to
look at $T_{\ell}$ as a vector space transformation and simply write
down the (matrix) elements then re-index them appropriately.
\end{enumerate}

\section{Corepresentations for $N>2$}
Before moving further it is necessary to realize the proper
generalizations of corepresentations on $SL_q(N)$ for $N>2.$  In the
case of $SL_q(2)$ and $SU_q(2)$ one is fortunate enough to have the
luxury of all relations being explicitly calculable (cf. [KS]
$\S$4.2.4). When one moves into the higher dimensional cases, one
quickly encounters a plethora of obstructions to calculating
everything explicitly.  Of course, it is possible to calculate
everything explicitly, but not in any concise manner.  The first
obstruction to note is that when one begins computing matrix
corepresentations for $SL_q(N)$ the index $\ell$ need not be
incremented by $1/2$ for each representation.  In fact notice that
even for $SL_q(3)$ if one considers $'T_{1/2}'$ the resulting matrix
is
\begin{equation}
T_{1/2} = \left(\begin{array}{ccc} u_{11} & u_{12} & u_{13} \\
u_{21} &u_{22} &u_{23} \\u_{31} &u_{32} &u_{33}
\end{array}\right) = \left(\begin{array}{ccc} t_{-1,-1} & t_{-1,0} & t_{-1,1} \\
t_{0,-1} & t_{0,0} & t_{0,1} \\t_{1,-1} & t_{1,0} &
t_{1,1}\end{array}\right).
\end{equation}

In particular $T_{1/2} = \{ t^{1/2}_{i,j}\}_{i,j=-1}^1.$ It turns
out that one needs to allow $\ell$ to increment appropriately.\\
\indent The appropriate increments of $\ell$ can be computed easily
using basic combinatorics.  When one needs to compute the matrix
corepresentation corresponding to $k$-homogeneous elements of
$\mathcal{O}(\C^N_q)$ i.e. $\mathcal{O}(\C^N_q)_k$ we have
$\left(\begin{array}{c} N+k-1 \\ k \end{array} \right)$ basis
elements.  In the specific case of $SL_q(3)$ one has
$\left(\begin{array}{c}k+2\\k\end{array}\right)
=\left(\begin{array}{c}k+2\\2\end{array}\right)$ elements, which one
will recognize easily as the familiar triangle numbers.\\
\indent How then should the re-indexing happen?  The task at hand
should not be so difficult if all of the indices were positive,
however in order to preserve as much useful information from the
$N=2$ case one should like to allow indices to run from $-\ell$ to
$\ell$ in unit increments.  For example one should like to have in
the case of $SL_q(3)$ acting on $\mathcal{O}(\C^3_q)_2$ to have a
$6\times 6$ corepresentation where indices should run from $-5/2$ to
$5/2$ in unit increments.  So one writes
\begin{equation}
T_{5/2} = \{ t_{i,j}^{5/2} \}_{i,j=-5/2}^{5/2}.
\end{equation}

One should also like to write down the correspondence
\begin{equation}T_{\ell}(e_s^{\ell}) = e_r^{\ell}\otimes t_{r,s}^{\ell}.
\end{equation}

Here the $e_s^{\ell}$ are the renormalized basis for
$\mathcal{O}(\C^N_q)_k$ with
\[
\ell =\frac{1}{2}(\left(\begin{array}{c}
N+k-1\\k\end{array}\right)-1).
\]

In order to give some consistency to (3.0.9) one requires a proper
indexing of $e_s^{\ell}.$  The explicit derivation of $s$ will be
given in the appendix.  For now, let it suffice to have a formula.
The given bases for $\mathcal{O}(\C^N_q)_k$ with renormalized
coefficients are
$$ \left(\begin{array}{c} k \\
(i_1,\dots,i_N)\end{array}\right)^{1/2}_{q^{-2}}\prod_{j=1}^N
x_j^{i_j}$$ and $$\{e_s^{\ell}\}_{s=-\ell}^{\ell}.$$

In order to equate these bases, the first important task is to give
an ordering to the first basis.  The proper ordering yields the
following map.

\begin{equation}
\left(\begin{array}{c} k \\
(i_1,\dots,i_N)\end{array}\right)^{1/2}_{q^{-2}}\prod_{j=1}^N
x_j^{i_j}  \mapsto e^{\ell}_s
\end{equation}

where
\begin{equation}
s = \sum_{r=1}^{N-1}\left(\begin{array}{c} \sum_{p=0}^{r-1}i_{N-p} + p \\
r \end{array}\right) - \frac{1}{2}(\left(\begin{array}{c} N+k-1 \\ k
\end{array}\right)-1)
\end{equation}

\begin{ex}
This formula looks treacherous, however it essentially gives an
inverse lexicographic ordering on the products of $x_i$ and lists
them in a bearably normal way.  For the case of
$\mathcal{O}(\C^3_q)_2$ one obtains the map
\begin{eqnarray}
(1+q^{-2}+q^{-4})^{1/2}x_1x_2 \mapsto e_{-3/2} & x_1^2 \mapsto e_{-5/2} \\
(1+q^{-2}+q^{-4})^{1/2}x_1x_3 \mapsto e_{-1/2} & x_2^2 \mapsto e_{1/2}\nonumber \\
(1+q^{-2}+q^{-4})^{1/2}x_2x_3 \mapsto e_{3/2} & x_3^2 \mapsto
e_{5/2}\nonumber
\end{eqnarray}
\end{ex}

\section{Appendix A: The Correspondence of Bases}
This section seeks only to show how the map between two bases of
$\mathcal{O}(\C^N_q)_k$ is obtained.\\
\indent To begin, note that one seeks an ordering on $
\left(\begin{array}{c} k \\
(i_1,\dots,i_N)\end{array}\right)^{1/2}_{q^{-2}}\prod_{j=1}^N
x_j^{i_j}$ and that $e_s^{\ell}$ is indexed by a single number so
that the ordering is easy. Consider the map
\begin{equation}
\left(\begin{array}{c} k \\
(i_1,\dots,i_N)\end{array}\right)^{1/2}_{q^{-2}}\prod_{j=1}^N
x_j^{i_j} \mapsto (i_1,\dots,i_N).
\end{equation}
One only needs to order the lists $(i_1,\dots,i_N)$ in some fashion.
The usual convention would be to take the simple lexicographic
ordering, but in this case the inverse ordering will be used so that
$(k,0,\dots,0)$ will correspond the leftmost matrix element.  With
this in mind consider the new map
\begin{equation}
(i_1,\dots,i_N)\mapsto r.
\end{equation}
The present goal is to compute $r$ and then to readjust $r$ into $s$
in a manner that allows $s\in \{
-\ell,-\ell+1,\dots,\ell-1,\ell\}.$\\
\indent The method to compute is is simply to set some parameters
and then read off how certain combinatorial moves affect the integer
ordering.  For example if $k=3$ and $N=3$ then the ordering will be
the following:
\begin{eqnarray*}
(3,0,0)\mapsto 0, &(2,1,0)\mapsto 1,&(2,0,1)\mapsto 2\\
(1,2,0)\mapsto 3, &(1,1,1)\mapsto 4, &(1,0,2)\mapsto 5,\\
(0,3,0)\mapsto 6, &(0,2,1)\mapsto 7, &(0,1,2)\mapsto 8,\\
(0,0,3)\mapsto 9
\end{eqnarray*}
Corresponding to 10 basis elements consistent with
$\left(\begin{array}{c}3+3-1\\3\end{array}\right) =10.$ One
conspicuous observation is that $(i_1,\dots,i_{N-1},i_N)\rightarrow
(i_1,\dots,i_{N-1} - 1,i_N + 1)$ corresponds to an increase of $1$
in the integer ordering. This is essentially the method of
observation used by the author to construct the ordering.  When one
looks back one step further to
$$(i_1,\dots,i_{N-2},i_{N-1},0)\rightarrow
(i_1,\dots,i_{N-2}-j,i_{N-1}+j,0)$$ one finds that this corresponds
to an increase in the ordering by $\sum_{p=0}^j p =
\left(\begin{array}{c}j+1\\2\end{array}\right).$  The rest of the
steps follow similarly. So that one finds at the stage
$$(\dots,i_{N-d},i_{N-d+1},\dots)\rightarrow
(\dots,i_{N-d}-j,i_{N-d+1}+j,\dots)$$ the increase is
$\left(\begin{array}{c}j+d\\d+1\end{array}\right)$ in the ordering.
Therefore a convenient way to write $(i_1,\dots,i_N)$ where $\sum
i_j = k$ is
\begin{equation}
(i_1,\dots,i_N) = (k-j_1,j_1-j_2,\dots,j_{N-1}-i_N,i_N).
\end{equation}
By which one can easily keep track of each combinatorial move.  In
short; where $(i_1,\dots,i_N)\mapsto r$ one has
\begin{equation}
r = i_n + \left(\begin{array}{c}j_{N-1}+1\\2\end{array}\right) +
\cdots + \left(\begin{array}{c} j_{N-1}+\cdots+j_1 +
(N-1)\\N-1\end{array}\right)
\end{equation}
and solving for each $j_d$ one obtains the formula
\begin{equation}
r = \sum_{r=1}^{N-1}\left(\begin{array}{c} \sum_{p=0}^{r-1}i_{N-p} + p \\
r \end{array}\right).
\end{equation}
The only thing left to do is to shift $r$ to $s$ so that $s$ is
symmetric about zero.  It has already been established though that
hen $\sum i_j = k$ there are $\left(\begin{array}{c} N+k-1 \\ k
\end{array}\right)$ basis elements.  One simply needs to subtract
one and divide by two to center this many numbers about zero thus
resulting in the rather horrendous formula as above
\begin{equation}
s = \sum_{r=1}^{N-1}\left(\begin{array}{c} \sum_{p=0}^{r-1}i_{N-p} + p \\
r \end{array}\right) - \frac{1}{2}(\left(\begin{array}{c} N+k-1 \\ k
\end{array}\right)-1)
\end{equation}

\end{document}